# Investigation of microstructural evolution of irradiation-induced defects in tungsten: an experimental-numerical approach

Salahudeen Mohamed[1]*, Qian Yuan[1], Dimitri Litvinov[1], Jie Gao[1]**, Ermile Gaganidze[1], Dmitry Terentyev[2], Hans-Christian Schneider[1], Jarir Aktaa[1]

[1]*Karlsruhe Institute of Technology (KIT), Institute for Applied Materials, Hermann-von-Helmholtz-Platz 1, 76344 Eggenstein-Leopoldshafen, Germany*

[2]*SCK • CEN, Institute of Nuclear Materials Science, Boeretang 200, B2400 Mol, Belgium*

*Corresponding author:*

Dr Salahudeen Mohamed

Karlsruhe Institute of Technology (KIT)

Institute for Applied Materials, Hermann-von-Helmholtz-Platz 1, 76344 Eggenstein-Leopoldshafen, Germany

Ph: +4915510920492

e-mail: Salahudeen.kunju@kit.edu

**Current affiliation:*- Fudan University, 220 Handan Rd, Yangpu District, Shanghai, China, 200437

**Abstract**

The hostile condition in a fusion tokomak reactor poses the main challenge in the development and design of in-vessel components such as divertor and breeding blanket due to fusion relevant irradiation conditions (14 MeV) and large thermal loads. The current work describes the employment of an integrated experimental-numerical approach to assess the microstructure evolution of dislocation loops and voids in tungsten proposed for fusion application. Cluster dynamics (CD) model is implemented and simulations are performed on the irradiated tungsten Disk shape Compact Tension (DCT) specimen used in the experimental test. TEM characterisation is performed on the DCT specimen irradiated at 400 °C and 600 °C with around 1 dpa, respectively. The dpa rate and cascade overlap rate from the experiments and SPECTRA-PKA code, respectively, are implemented in the CD model. Based on the comparison between experimental and computational results, the dose and temperature dependence of irradiation-induced defects (dislocation loops, voids, c15 clusters) are clearly observed. Trap mediated diffusion is studied and the impact of cascades with the pre-existing defects is analysed through full cascade overlap mode and the consequent influence on the defect concentration is evaluated. The exchange of self-interstitial atoms (SIAs) and the change in the size of loops through reaction between ½<111> and <100> loops are studied in detail by means of the transfer rate of the SIAs.

# 1  Introduction

The future tokamak fusion reactors, International thermo-nuclear experimental reactor (ITER) and DEMOnstration tokamak (DEMO), are installed with in-vessel components which are subjected to high fluence of 14 MeV neutrons and high heat flux loads [1,2]. Tungsten (W) is employed as armour material in the divertor of ITER and DEMO due to its low sputtering rate, high melting temperature and high thermal conductivity [3,4]. The hostile fusion relevant irradiation environment in the tokamak reactor during operation leads to the formation of irradiation-induced defects in W such as dislocation loops, voids and transmutation products due to Primary Knock-on Atoms (PKAs) [5,6]. Consequently, W material undergoes degradation due to irradiation-induced cascade damage and transmutation which alters the thermo-physical and mechanical properties at different stages of its lifecycle during operation [7,8]. In particular, the irradiation-induced defect formation in W leads to the various issues such as hardening, swelling and embrittlement. Moreover, these irradiation-induced defects contribute to the increase in ductile-to-brittle transition temperature (DBTT) in W [9]. Due to these factors W based components in tokamak are prone to failure, which consequently reduces its lifespan [10]. Therefore, efforts are required to gain some insights into the formation and evolution of the irradiation-induced dislocation loops and voids and utilize this information for the development of W based components for fusion relevant applications.

Several studies are carried out by many authors using experimental analysis to investigate the underlying mechanisms on the microstructural evolution of irradiation-induced defect features in W [11–15]. In the study carried out by Li et al., irradiation experiments using 30 keV He$^+$ were conducted on W samples to understand the evolution of dislocation loops depending on irradiation temperature, dose and thickness of the sample using *in-situ* TEM [14]. The study demonstrated that the pre-existing dislocation lines have a significant influence on the Burgers vector, density and size of the dislocation loops. Moreover, the role of the surface of thin foil acting as sinks influences the overall density distribution of the dislocation loops. Using TEM analysis, Klimenkov et al., investigated the neutron irradiated W about 1 dpa at various irradiation temperatures and observed irradiation-induced dislocation loops and voids [16]. In the experimental work carried out by Chauhan et al. on the neutron irradiated W about 1 dpa, coarsening of voids is observed when increasing post-irradiation annealing temperature and time [15]. Dislocation loops and voids of size 10 nm and 10 nm – 65 nm, respectively, are examined by TEM by Duerrschnabel et al. in W material irradiated at 800 °C and 1.25 dpa. The dislocation loops populations of ½<111> and <100> are reported in various experimental studies [17,18]. ½<111> loops only are observed in W irradiated with 15-85 keV He$^+$ to 3 dpa at 500 °C in [19–21]. However, in some experimental studies, coexistence of ½<111> and <100> loops are noticed. In body centered cubic (bcc) materials, the <100> loops are almost immobile when compared to the highly mobile ½<111> loops. The less mobility of <100> loops contributes to the irradiation hardening and embrittlement in W [14].

The irradiation of the in-vessel materials is currently conducted in surrogate fusion reactor facilities since the actual fusion environment material testing facilities are not yet available. Moreover, employing only experimental testing to sort and down select the candidate material of in-vessel components is expensive. In order to accelerate the fusion relevant material development, predictive models are useful which can emulate fusion irradiation conditions. It could aid in the reduction of costs related to conducting irradiation experiments. In particular, the irradiation damage in nuclear materials is multi-scale in nature and there is a need to understand the long-term kinetic evolution of defects. In this aspect, mesoscale models like cluster dynamics (CD) and kinetic Monte Carlo (KMC) methods are usually employed to analyse the evolution of irradiation-induced defect features in W [22]. In the case of KMC approach, due to the constraints in the computational domain, which can accommodate only a defect of few nanometers while for CD method, there are no constraints to incorporate the size of the defect clusters and it can handle defects with larger cluster sizes. However, there is a need for higher computational resources in CD to perform simulation for larger defect clusters. To overcome these difficulties, various mathematical schemes are developed to reduce the cost associated with the CD method. In the work of Zhao et al., CD methods are employed to study the helium/hydrogen retention in W [23]. In the numerical analysis conducted by Li et al. [24] and Krasheninnikov et al. [25], the nucleation and diffusion of helium atoms in W are analysed in detail using CD model. Using spatial dependant CD method, Shah et al. [26] have analysed the formation of the helium vacancy clusters in W monoblock and found that due to emission of helium from the surface, lower retention of helium clusters is observed. A similar study is conducted by Faney et al. [27] using spatial dependant CD method to analyse the tendrils formation in W exposed to helium.

So far, based on the previous computational modelling works using CD method in W, there is no study, which investigated the long-term evolution of dislocation loops and voids. In order to address this aspect, an in-house CD method has been developed, which has the flexibility to integrate defects features such as dislocation loops, voids, c15 clusters and impurities as a function of irradiation dose and temperature. The present study aims at investigation of the microstructural evolution of irradiation-induced defects using an integrated experimental-numerical approach. Moreover, interaction of loops of different Burgers vectors, ½<111> and <100>, with voids, impurities and transmuted elements and the resulting influence on the overall evolution of defect populations is studied. Dislocation loop growth due to the reaction between ½<111> and <100> loops and trap mediated 1-D loop migration are also considered.

The paper is organised as follows: The next section describes the CD model; Section 3 reports the experimental results in comparison to those obtained using the CD based numerical simulations.

## 2 Model description

### 2.1 CD model

Cluster dynamics (CD) method is a mathematical model based on a rate theory approach to mimic and study the kinetics of microstructural evolution of defects in materials [24,28,29]. It provides insights on the mechanisms involved in the radiation damage, which can be employed to analyse nucleation of dislocation loops, void formation, precipitation and several other irradiation-induced phenomena like swelling [30]. In CD method/simulation, the microstructural evolution of the system is described by means of a system of differential equations on the concentration of the cluster. It computes evolution of the cluster size and concentration through space and time based on the kinetic reaction of irradiation-induced defects considered in the system. The evolution of cluster concentration of irradiation-induced defects is represented by the following equation:

$$\frac{dC(y)}{dt} = G_{V,I} - L^y_{sink} - \sum_x k^x_y C(y) C(x) - \sum_x g^x_y C(y) + \sum_x k^x_{y'} C(y') C(x) + \sum_x g^x_{y''} C(y'') \quad (1)$$

where $G_{V,I}$ denotes defect production rate. $C(x)$ and $C(y)$ represent concentration of defects of types $x$ and $y$, respectively, $L^y_{sink}$ signifies the loss rate of the defects through sinks (surface, dislocation lines, grain boundaries). The interaction among the irradiation-induced defects are defined in terms of reaction coefficients ($k_y^x$, $g_y^x$, $k_{y'}^x$ and $g_{y''}^x$). $k_y^x$ describes the rate of defect $y$ absorbing defect $x$ while $g_y^x$ represents rate of defect $y$ emitting mobile defect $x$. $k_{y'}^x$ denotes the generation of defects $y$ due to the reaction between defects $y'$ and $x$. $g_{y''}^x$ depicts the emission of defect x by defect $y''$.

As a result of neutron irradiation, the cascade events occur in the material. To reproduce these damage events in CD model, the defect production rate has to be employed. Based on the work of Sand et al., the defect production term is expressed by a power law which relates the cluster size, $n$, and defect production rate, $G_{V,I}$ [31,32] as shown in equation (2):

$$G_{V,I} = \frac{A_{V,I}}{n^{S_{V,I}}} \quad (2)$$

$V$ is vacancy, $I$ is self interstitial, $A_{V,I}$ is the pre-factor related to total production rate, $G_{dpa}$. $S_{V,I}$ denotes the scaling component [26,33,34]. The pre-factor, $A_{V,I}$, of the power law is calculated based on the total defect production rate, $G_{dpa}$, and it reads:

$$A_{V,I} = \frac{G_{dpa}}{\sum_{n=1}^{N_{max}} n^{1-S_{V,I}}} \quad (3)$$

The parameters $S_V$ and $S_I$ are 1.63 and 2.20, respectively, obtained from the work of Sand et al. [31,32]. $N_{max}$ for interstitials and vacancies are taken as 15 and 20, respectively [35]. To incorporate the reduction of the defect production rate with time, the defect production rate ($G_{dpa}$) decreases linearly based on the defect fraction ($f_{V,I}$) in the materials [26,34,36]. Therefore, the final form of equation is written as:

$$G_{V,I} = \frac{A_{V,I}}{n^{S_{V,I}}}(1-f_{V,I}) \quad (4)$$

However, for the current study, $G_{dpa}$ is maintained uniform for the experimental analysis. Therefore, $f_{V,I}$ is zero for the CD simulation. The irradiation-induced defects in the CD simulation domain consist of self-interstitial atoms (SIAs), self-interstitial dislocation loop types of different Burgers vectors (½<111>, <100>), c15 clusters, vacancy and voids (vacancy clusters) [37–39].

SIA clusters comprising less or equal to 3 SIAs and vacancy clusters containing vacancies below 4 (<=4) have 3-D motion. SIA clusters with size greater than 3 are considered as ½<111> loops, <100> loops, c15 clusters, which perform 1-D migration. Voids containing more than 4 vacancies and c15 clusters containing SIAs are treated as immobile defects. Vacancy type dislocation loops have not been taken into account in the current study [37].

*Table 1. CD model mobility parameters for interstitial and vacancy defects.*

| Cluster | Pre-exponential factor (m$^2$/s) | Migration energy (eV) |
|---|---|---|
| $I_1$ (3-D) | $D_0(I_1)$ : $9.981 \times 10^{-11}$ [40] | 0.165 [40] |
| $I_2$ (3-D) | $8.648 \times 10^{-10}$ [40] | 0.222 [40] |
| $I_3$ (3-D) | $3.47 \times 10^{-11}$ [40] | 0.203 [40] |
| $I_4$-$I_\infty$ (1-D) | $D_0(I_1)/n^{0.5}$, n is the number of SIAs in ½<111> loop [41] | 0.1 [40] |
| $V_1$ (3-D) | $D_0(V_1)$ : $177 \times 10^{-8}$ [41] | 1.66 [41] |
| $V_2$ (3-D) | $D_0(V_1)/2$ [37] | 1.66 [41] |
| $V_3$ (3-D) | $D_0(V_1)/3$ [37] | 1.66 [41] |
| $V_4$ (3-D) | $D_0(V_1)/4$ [37] | 1.66 [41] |
| $V_5$-$V_\infty$ (1-D) | 0 | - |

The reaction for 3-D configurational structures like c15 clusters and 3-D migrating defects ($I_1$-$I_3$, $V_1$-$V_4$) are modelled using absorption coefficients, which are calculated using the following equation:

$$k_y^x = 4\pi(D_y + D_x)(r_y + r_x) \tag{5}$$

$D$ is diffusion coefficient and $r$ is the reaction radius.

For the reaction between 3-D migrating defects ($x$) and loops ($y$) involving reactants with non-spherical volumes, reaction rate is evaluated based on equation (6) [42]:

$$k_y^x = (D_y + D_x)((1 - \alpha_y^x) z_{y,x}^L + \alpha_y^x z_{y,x}^V) \tag{6}$$

where,

$\alpha_y^x = \left(1 + \left(\frac{r_y}{3(r_x+r_d)}\right)^2\right)^{-1}$, $r_d$ = dislocation capture radius

$z_{y,x}^V = 4\pi(r_y + r_x + r_d)$

$z_{y,x}^L = \frac{4\pi^2 r_y}{\log(1+8r_y/(r_x+r_d))}$

The interaction between 1-D migrating defects is modelled using the reaction rate [42].

$$k_y^x = \frac{2 D_y \sigma_{y,x}}{\lambda_y}$$

$\lambda_y$ represent mean free path of the 1-D migrating defect calculated by $\lambda_y^{-1} = \sum_\omega \sigma_{y,\omega} C(\omega)$, where $\omega$ denotes the reaction of all defects and sinks with 1-D migrating defects. $\sigma_{y,x}$ denotes the cross-section of the reactants [42].

The emission coefficient $(g_y^x)$ is calculated via the binding energy, $E_b^x$, and absorption coefficient, $k_{y-x}^x$, and it reads:

$$g_y^x = \frac{1}{\Omega} k_{y-x}^x \exp\left(\frac{-E_b^x(y)}{k_B T}\right) \tag{7}$$

where $k_B$ is Boltzmann's constant.

In order to obtain the binding energy of the single SIA and vacancy with the SIA and vacancy clusters, equation (8) - equation (11) are employed. $E_f(n)$ is the formation energy of SIA clusters or voids, which consists of $n$ SIAs or voids as shown in the following equations.

*Vacancy cluster*:

$$E_b^{SIA}(n) = E_f^{SIA} + E_f(n+1) - E_f(n) \tag{8}$$

*SIA cluster*:

$$E_b^{SIA}(n) = E_f^{SIA} + E_f(n-1) - E_f(n) \tag{9}$$

*Vacancy cluster*:

$$E_b^V(n) = E_f^V + E_f(n-1) - E_f(n) \tag{10}$$

*SIA cluster*:

$$E_b^V(n) = E_f^V + E_f(n+1) - E_f(n) \tag{11}$$

Where $E_f^V = 3.80$ eV [43], $E_f^{SIA} = 9.46$ eV [44].

The formation energy of loops of ½<111> and <100> types are taken from the ab-initio study of Alexander et al.[45] using the analytical expression:

$$E_f(n) = a_3 + a_2\sqrt{n} + a_1\ln(n)\sqrt{n} \tag{12}$$

The coefficients $a_0$, $a_1$, $a_2$ for ½<111> and <100> loops are:

*½<111> loop:*
$a_1$= 3.92996
$a_2$= 7.92419
$a_3$= 6.20090

*<100> loop:*
$a_1$= 4.84883
$a_2$= 13.6984
$a_3$= -8.2584

In the case of voids, formation energy is evaluated as:

$$E_f^V(n) = 4\pi(r_V)^2\gamma, \tag{13}$$

where $\frac{4}{3}\pi(r_V)^3 = n\frac{a_0^3}{2}$, $a_0$ is lattice parameter, $\gamma$ is surface energy.

Grain boundaries, dislocation lines and free surfaces are employed as sinks in the current model for 1-D and 3-D migrating defects [27]. The sink strengths for the defects implemented in the current study are described as follows:

*Sink of Dislocation lines* [46]:

$$L_{dl,1D} = 2(\pi(r_d+r_{V,I})\rho_d)^2 \text{ (1-D migrating defect)} \tag{14}$$

$$L_{dl,3D} = \frac{2\pi\rho_d(1-\rho^2)}{\ln(\frac{1}{\rho})-0.75+0.25\rho^2(4-\rho^2)} \text{ (3-D migrating defect)} \tag{15}$$

where $\rho = (r_d+r_{V,I})\sqrt{\pi\rho_d}$

*Sink of Grain boundary* [46]:

$$L_{GB,1D} = \frac{15}{R_{GB}^2} \text{ (1-D migrating defect)}$$

$$L_{GB,3D} = \frac{14.4}{R_{GB}^2} \text{ (3-D migrating defect)}$$

*Sink of free surface* [47]:

$$L_{surface,1D} = \frac{8\cos^2\varphi}{l^2}\text{(1-D migrating defect)}$$

$$L_{surface,3D} = \frac{2}{l^2}\text{(3-D migrating defect)}$$

where $l$ is the thickness of the sample, $\varphi$ is the angle of direction of loop migration and surface normal.

In experimental studies, it has been shown that the irradiation-induced defects like dislocation loops are trapped by interstitial or substitutional atoms which thereby imped its motion [16]. Within the trap sites, loops perform random walks in addition to the movement among lattice sites and in order to take into account this phenomena, an additional diffusivity term, $D_i^{dt}$, is incorporated based on equation (16) [30].

$$D_i^{dt} = v_i^{dt}\frac{\lambda_i^2}{2N} \tag{16}$$

where $\lambda$ is the hop length to account for the jump between trapping sites which is assumed to be one dimensional (N=1), $v_i^{dt}$ is the activation frequency. The hop length, $\lambda_i$ is proportional to the trap density, $\rho_t$, and $r_t$ is the length between trap and mobile dislocation loop [30] as shown in equation 17.

$$\lambda_i = (4\pi r_i r_t \rho_t)^{-1} \tag{17}$$

½<111> loops are treated as trapped and freely migrating loops in the present study. The free ½<111> loops are modelled in CD in such a way that they can be trapped by the trapping elements in the course of motion.

*Table 2. CD model input parameters for W.*

| Parameters | Value |
|---|---|
| Lattice parameter, $a_0$ | 0.316 nm [26] |
| Atomic volume | $a_0^3/2$ [24] |
| Burgers vector, b | 0.274 nm [24] |
| Dislocation density, $\rho_d$ | 4-6 x $10^{12}$ m$^{-2}$ [26] |
| Grain size, d | 3 µm [26] |
| Dislocation capture radius, $r_d$ | 0.65 nm [26] |
| Dislocation bias, Z | 1.2 [26] |
| Trap radius, $r_t$ | 1.5b [42] |

In the study of Liu et al. [39], during cascade, majority of the dislocation loops formed in W is of interstitial type. In particular, different groups of interstitial clusters are categorized: i) c15 clusters ii) ½<111> loops iii) <100> loops iv) Mixed structures comprising different Burgers vector v) 'None' structures. At all PKA energies, most of the interstitial clusters observed are ½<111> loop type while at high energy cascades, a few <100> loop type and mixed structures are noticed. Furthermore, the probability of c15 cluster formation increases with an increase of the PKA energy. In another study conducted by Liu et al. [38], a detailed MD analyses were carried out to analyse the formation and stability of c15 clusters. In their study, it has been shown that c15 clusters are immobile and highly stable once nucleated in W and can collapse at high temperatures above 1500°C to transform to ½<111> loop type [38]. Moreover, in the event of collapse of c15 clusters, before the conversion of c15 cluster into ½<111> loops at higher temperatures, the c15 clusters can transform into <100> loops at lower temperatures [38]. This is due to the smaller energy difference in the formation energy values between ½<111> and <100> loops. Based on the above study, the current work has employed the formation of interstitial clusters, ½<111> and <100> loops, SIA clusters of 2 and 3 SIAs and c15 clusters during irradiation cascades. In particular, based on the work of Liu et al [38,39], the study has assumed the formation of 5% c15 clusters, 1% <100> loops and the remaining as <111> dislocation loops and clusters comprising 2-3 SIAs. The dislocation loops (½<111>, <100>) and c15 clusters have greater than 4 SIAs. However, mixed and 'None' structures are not taken into account for the irradiation cascade due to the insufficient data. In order to incorporate their role in the evolution of the defects, the percentage fraction of these interstitial structures are included in the dislocation loops (½<111>, <100>) and c15 clusters based on the similar approach employed in the study of Gao et al (Ref. [37], figure 4)

In the current study, the loop evolution in W by absorption and coalescence through interaction between the loops of same and different Burgers vectors is considered as shown in Figure 1. The following loop interactions are modelled in this study based on the former studies. The dislocation loops of interstitial type increase its size by absorbing SIAs or smaller mobile loops since the SIA have higher strain field than vacancies [13,48]. Moreover, the growth of highly mobile ½<111> loops occur by absorbing or coalescing with smaller <100> loops. Similar mechanisms are considered for the <100> loops in which larger size of <100> loops can grow by absorbing smaller ½<111> loops. In fact, the growth of ½<111> and <100> loops are mainly by transfer of SIAs between them. The reaction between different ½<111> loops variants containing same number of SIAs to obtain <100> loops is rare. However, the current study has incorporated this scenario and assigned a probability of ∼3/8 (3c/(8c-2)), c is the concentration of each variant [37]. The remaining 5/8 between ½<111> loops of different variants and identical size yields ½<111> loops of larger sizes [37] .

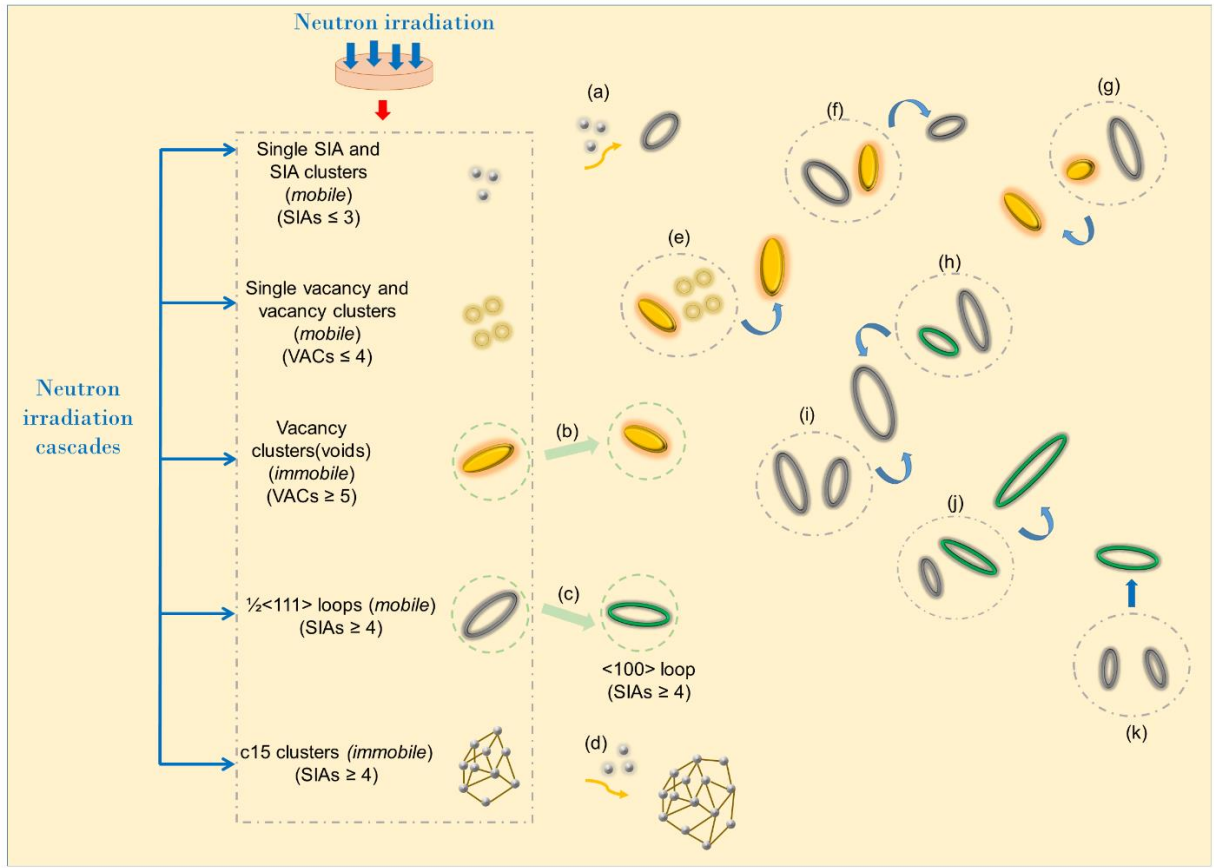

*Figure 1. Illustration of defect kinetics involved in the evolution of irradiation-induced dislocation loops and voids. (a) depicts the ½<111> loop nucleation due to SIAs migration, (e) depicts the nucleation of voids via migrating mobile vacancies, (b)-(c) illustrates the cascade overlap with pre-existing voids and dislocation loops, which were formed during processes (a) and (e), (d) represent the nucleation of c15 clusters due to SIAs migration, (f)-(k) depicts the voids or loop coalescence caused by 1-d loop migration.*

The presence of interstitial impurity atoms and substitutional transmuted atoms in W plays an important role in the evolution of irradiation-induced void and loops. In the case of pure W grade employed in the present study, despite having 99.7% purity, interstitial impurities such as C, H, He, N and O are observed with concentrations reported in the study of Bakaev et al [49]. Moreover, transmuted elements such as Rhenium (Re) atoms are observed in the neutron irradiated W sample employed for the irradiation reported in the previous and current studies [15,50]. The impurities bind with vacancies and interstitial clusters and impedes the movements of them. Therefore, the impurities act as traps for the irradiation-induced mobile defects, which influence the overall microstructure of the material. It is thus important to incorporate interstitial impurity atoms and substitutional transmuted atoms in the current CD model. In the current study, trap elements of concentration of 229 appm are employed as interstitial impurity atoms [49]. In addition, transmutation based atoms of density of $1.2 \times 10^{14}$ cm$^{-3}$ [15] are also modelled as trap elements. Moreover, the CD model considers only the interaction of ½<111> loops and mobile vacancies with the traps.

*2.2    Numerical algorithm*

The model employs the Jacobian matrix of the ODE system using the analytical form, $J_{ij} = \partial(\frac{\partial C_i}{\partial t})/\partial C_j$, Where $C_i$ and $C_j$ are the density of clusters. The solution of the Jacobian matrix is evaluated by means of python wrappers of the odespy package [51] based on the odepack algorithms [52]. In order to model large defect clusters, a higher number of differential equations has to be coupled to calculate the solution. To speed up the computational process and time, grouping method is employed, which implements two equation per cluster group to achieve conservation of mass and the cluster density [53].

## 3    Results

*3.1    Experimental parameters implemented for CD model*

For the experimental analysis, ITER graded stress-relieved IGP-W bar, dimensions $36 \times 36 \times 480$ mm$^3$ with 99.7 wt% purity produced by PLANSEE SE, Austria, is implemented [50]. Disk shape Compact Tension (DCT) W specimens with 4 mm thickness are fabricated using electrical discharge machining (EDM) transverse to forging direction [50]. Irradiation experiments are performed on DCT specimens placed in the fuel channel of Material Test High Flux BR2 reactor of SCK CEN in Mol. Experiments are performed for irradiation doses of ~ 1 dpa with neutron fluence between 0.1 MeV and 1 MeV at irradiation temperatures of 400 °C and 600 °C followed by TEM characterization of irradiation-induced defect features. The relevant experimental conditions are depicted in the table 3.

*Table 3. Irradiation conditions in experimental analysis.*

| $T_{irr}(^oC)$ | Irradiation time (days) | Flux (n/cm$^2$/s) (E>0.1 MeV) | Flux (n/cm$^2$/s) (E>1 MeV) | Flux (n/cm$^2$/s) (0.1<E<1 MeV) |
|---|---|---|---|---|
| 400 | 186 | $1.9 \times 10^{14}$ | $8.9 \times 10^{13}$ | $1.0 \times 10^{14}$ |
| 600 | 186 | $2.3 \times 10^{14}$ | $1.1 \times 10^{14}$ | $1.2 \times 10^{14}$ |

*3.2    TEM characterization of irradiation defects*

For microstructural investigations, a FEI Scios focused-ion beam (FIB) scanning electron microscope (SEM) was utilized for preparing TEM lamellae from the undeformed region of the disk shape compact tension (DCT) specimens. Dislocation loops and voids were quantitatively characterized by extensive transmission electron microscopy (TEM) investigation. For more details on sample preparation and TEM investigation, please see Ref. [15]. In addition, analytical TEM was used to examine the spatial distribution of the transmutation products. Figure 2(a) shows an overview of dislocation loop microstructure in IGP W after irradiation to 1 dpa at 400 °C. For loops quantitative analysis, several images were acquired and analysed under KBF conditions. For instance, a representative KBF micrograph (Figure 2 (b)) taken under g = {110} diffraction condition presents clear dislocation loop features with mainly elliptical shape where the major axis of the loop is taken as a measure of its size. Here, g = {110} and g = {002} diffraction vector was analysed, resulting in the loop densities of 6.6 $\times 10^{15}$ cm$^{-3}$ and 7.2 $\times 10^{15}$ cm$^{-3}$ and mean sizes of 3.6 nm and 4.8 nm, respectively. By applying the loop invisibility criteria and the statistical method for Burgers vector determination, the visible loops were found to be of both ½ <111> and <100> types, with 47% and 53% fraction, respectively. Owning to the resolution limit under KBF condition, loop size less than 4 nm is excluded in the statistical analysis of <100> loop fraction. On the basis of this KBF analysis, the total loop density and loop mean size were determined to be ~1.1 $\times 10^{16}$ cm$^{-3}$ and 4.1 nm, respectively. The size distribution of the dislocation loops under KBF condition with g = {110} diffraction vector is presented in Figure 5 (b). Figure 2 (c) presents the overall microstructure of the as-irradiated (1.06 dpa, 600 °C) IGP W acquired via weak-beam dark-field (WBDF) technique under g = {002} diffraction condition manifesting a high density of nm-size black dots and dislocation loops. The latter are better resolved in Figure 2 (d) at higher magnification. Typical loop features with circular/elliptical shape and coffee bean contrast are identified. For loops statistical analysis, total dislocation loop density and mean size are found to be ~$5.9\times10^{16}$ cm$^{-3}$ and 2.3 nm, respectively. Furthermore, as most defects (≤2 nm) show black dot features, which can be hardly resolved in TEM even via WBDF technique, only defects with characteristic features (Figure 2(d)) were counted as dislocation loops. Moreover, the <100> loop fraction is found to be about 21% which is calculated via the statistical method for Burgers vector determination based on the invisibility criteria (g·b=0). The size distribution of dislocation loops determined for g = {002} diffraction vector is presented in Figure 5 (c).

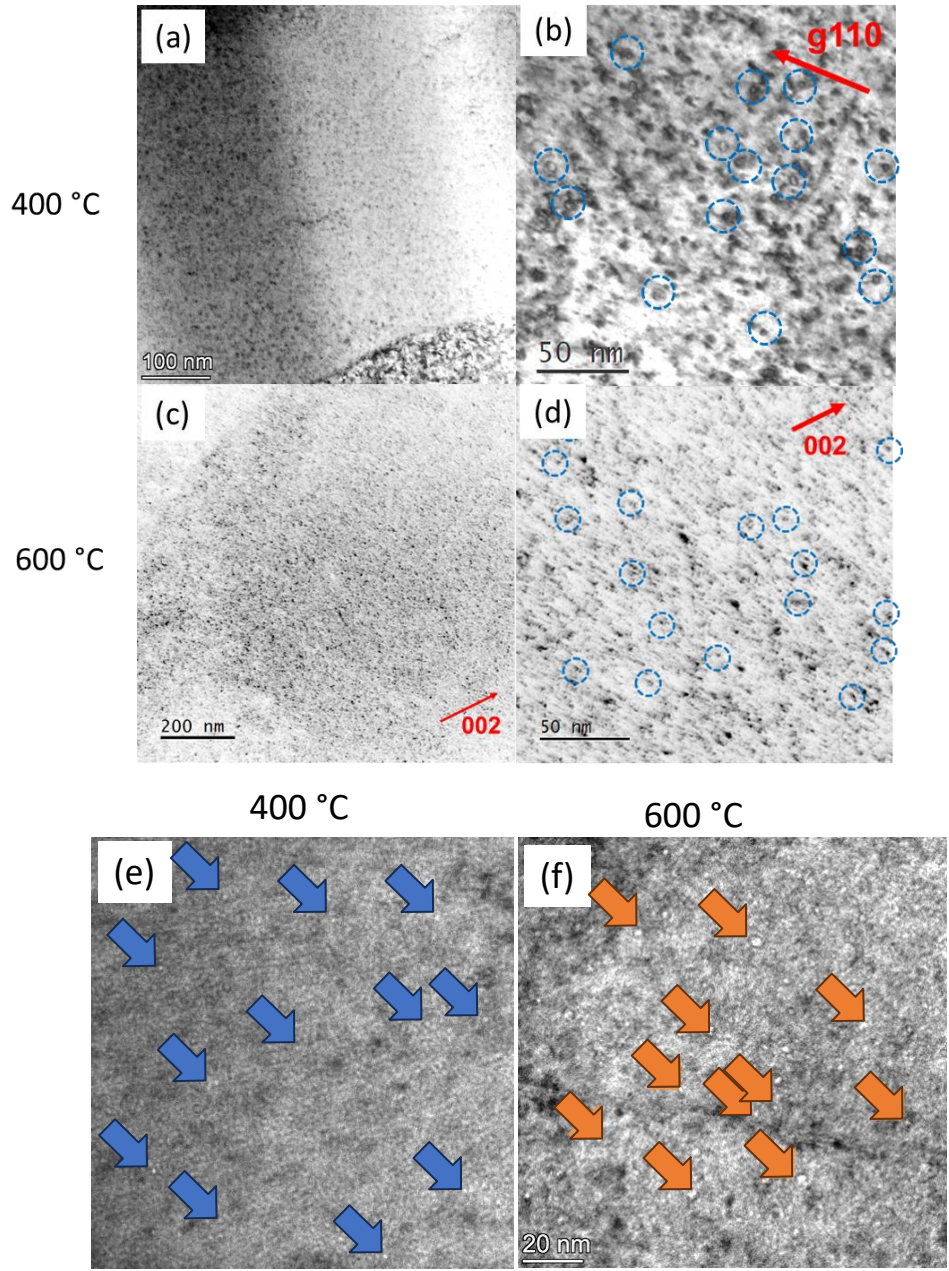

*Figure 2. Dislocation loop microstructure in irradiated IGP W (a) Inverted-contrast HAADF-STEM overview micrograph (1.02 dpa, 400 °C); (b) KBF micrograph taken under g = {110} diffraction condition; (c) Inverted-contrast WBDF micrograph taken under g = {002} diffraction condition showing an overview of as-irradiated (1.06 dpa, 600 °C) microstructure of W, which manifests a high density of black dots and dislocation loops; (d) the same as in (c) at a higher magnification for better resolution of dislocation loops with near circular, elliptical and coffee bean contrast (encircled); Typical BF-micrographs presenting void distribution in irradiated IGP-W at (e) 1.02 dpa, 400 °C (f) 1.06 dpa, 600 °C. Voids with sizes less than 1 nm were not counted in the statistics.*

The identification of void and its statistics is obtained by the traditional through-focal series technique where the contrast of void changes from a white dot with a black fringe in an under-focused image to a dark dot with a white Fresnel fringe in the over-focused image. Similar to the high density of dislocation loops, majority of nano-sized voids are found to be uniformly distributed in the as-irradiated sample, as presented in the under-focused BF-TEM micrograph in Figure 2 (e) for 400 °C. The mean voids diameter (fringe edge-to edge distance) is estimated to be 1 nm and the density is determined to be $1.9 \times 10^{17}$ cm$^{-3}$. Figure 7 shows the voids size distribution in which most voids are less than 2 nm. In addition, it is noteworthy that many voids with sizes less than 1 nm are hardly determined via through-focal series technique. Therefore, such ultra-small voids are excluded in the statistical analysis in this work.

For 600 °C, in addition to a few dislocation loops/black dots, several marked voids are identified near the grain boundary (Figure 2 (f)). This indicates that there is no clear void denuded zone near the grain boundaries in the as-irradiated condition. The mean void diameter (fringe edge-to edge distance) is estimated to be 1.5 nm and the density is determined to be $\sim 7.1 \times 10^{16}$ cm$^{-3}$. The size distribution of voids is presented in Figure 7.

As the content of transmutation products Re (Os) in the as-irradiated condition is expected to be at the level of 2 at% (0.2 at%), STEM-EDX analysis was conducted to visualize their distribution. No evidence of rhenium (Re) or osmium (Os) decoration or segregation at grain boundaries was observed at 400 °C, as depicted in Figure 3. Additionally, there was no evidence of precipitation within the grains. Notably, the as-irradiated sample exhibited approximately 2.4 at% Re and 1.5 at% Os, as confirmed by elemental mapping.

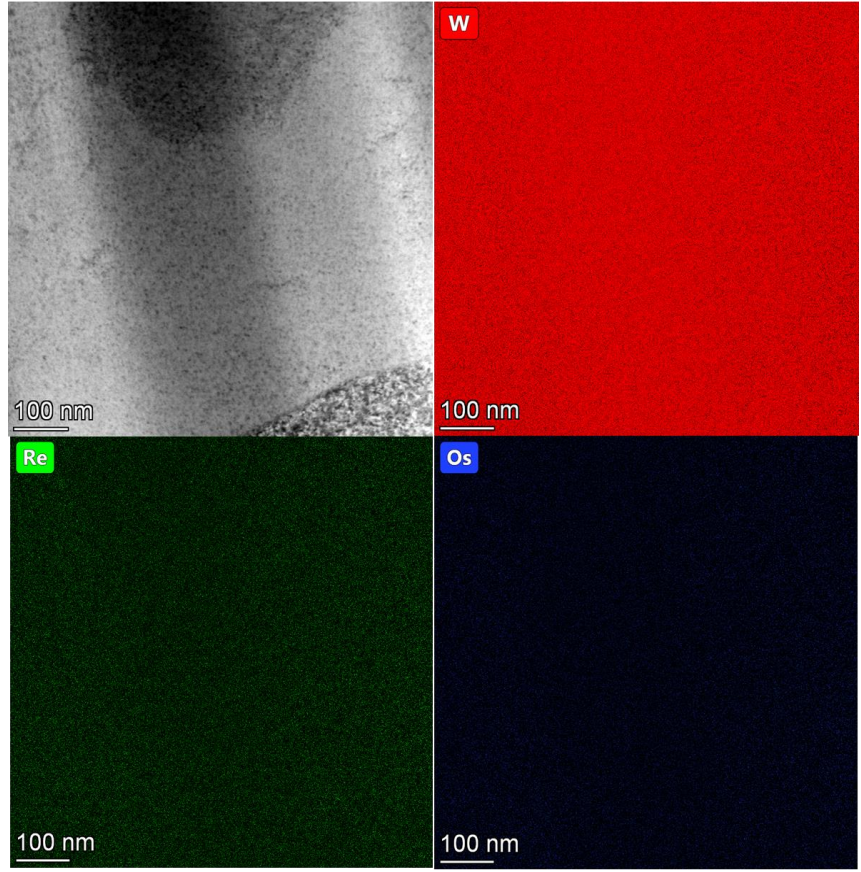

*Figure 3. Representative STEM-HAADF micrographs with inverted contrast, for irradiated (1.02 dpa, 400 °C) IGP-W along with corresponding STEM-EDX spectrum images providing elemental mapping, with W highlighted in red and transmutation-induced rhenium (Re) in green, along with osmium (Os) in blue.*

Elemental mapping (Figure 4) shows neither Re/Os decoration/segregation nor precipitation at grain boundaries and/or in the grains of the sample after irradiation at 600 °C. In addition, elemental quantification analysis shows that about 1.9 at% Re and 0.6 at% Os are detected in the selected grain (see marked frame in Figure 4), which indicates that Re , Os and their TEM-invisible clusters are uniformly distributed in the sample.

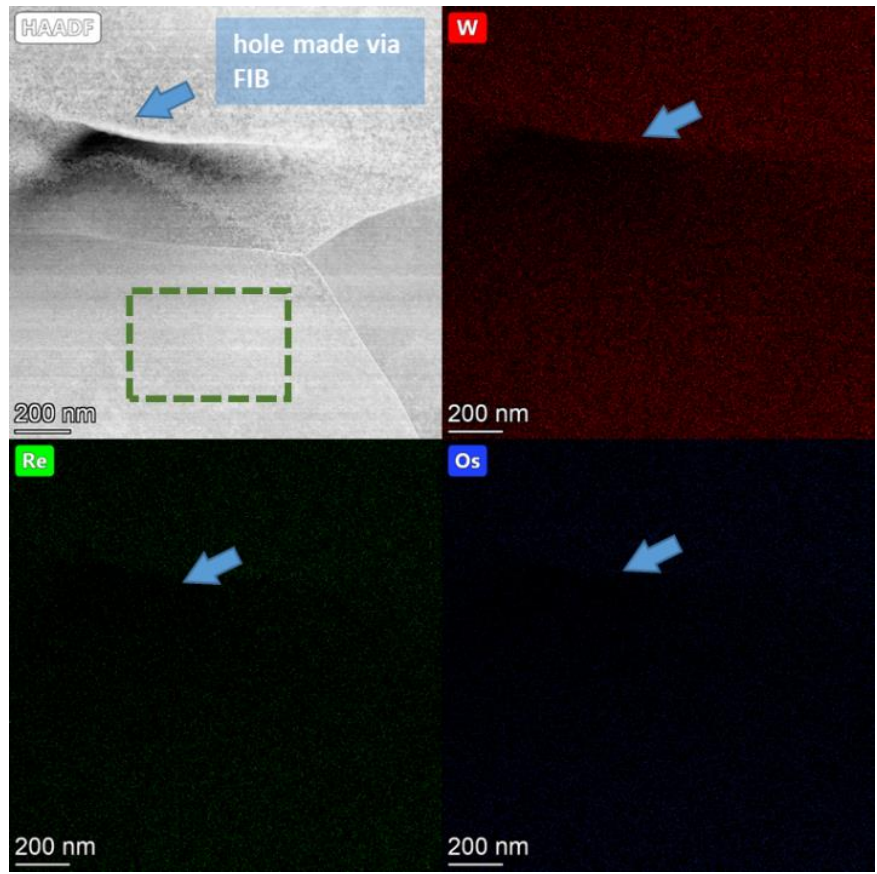

*Figure 4. Representative STEM-HAADF micrograph of irradiated IGP W (1.06 dpa, 600 °C) showing grain boundaries. STEM-EDX spectrum images show the corresponding elemental mapping of W (red) and transmutation induced Re (green) and Os (blue).*

*3.3    Cascade parameters implemented for CD model*

Irradiation damage in nuclear materials caused by the collision cascades is reported in several studies [54,55]. During continuous irradiation on the material, the collision cascades overlap with the pre-existing defects that can actually change the morphology and the size of the defect clusters in the irradiated material [54]. In fact, there is a reduction in the defects formed during cascade overlap in an irradiated material with respect to the amount of defects produced in a pristine material during the initial cascades, which is corroborated in the various studies [54,55]. Moreover, the cascade overlaps are categorised into partial overlap and full overlap [54]. In the case of partial cascade overlap, there is a drastic decrease in the defect cluster size while for full cascade overlap the average decrease in the defect cluster size is lower. The current study has incorporated the influence of the full cascade overlap on the pre-existing defect distribution on the W material. In order to realize cascade overlap event, transformation of clusters (½<111>and <100>) of SIA (<50) type is implemented in the CD model by means of transformation coefficients as implemented in the following governing equation for ½<111> and <100> loops (equation 18). The transformation coefficients represent the probability of transformation of the ½<111>and <100> loop into <100> or ½<111> loop types. The information regarding the transformation coefficients are taken from the study of Byggmastar et al. [54]. In the event of cascade overlap, there are also changes in the defect cluster morphology with the formation of mainly ½<111>, <100>, mixed and 'None' clusters. In fact, mixed clusters include dislocation loops of different Burgers vector (½<111>, <100>). The presence of mixed clusters are observed in the collision cascade simulations [54]. However, these mixed clusters are unstable and can transform into other loop types after annealing. In the case of 'None' structures, there are no information regarding its stability during a cascade overlap event in W. Due to scarce data on the mixed and 'None' clusters regarding its morphology and size, the current study has not employed its role in the cascade overlap event. Moreover, the cascade overlap event on the c15 clusters, and its transformation are not considered since there are insufficient information on the transformation coefficients.

In the current work, the influence of full cascade overlap is taken into account to understand the microstructural evolution of loops and voids in W. In order to incorporate the full cascade overlap event, cascade overlap rate ($F_{ol}$) is included in the CD model which considers the probability of full cascade overlap event per unit time based on equation 18 [37]:

$$\frac{dC_{111}(n)}{dt} = \epsilon + f_{111to100}(n)F_{ol}(n)C_{100}(n) - f_{100to111}(n)F_{ol}(n)C_{111}(n) \qquad (18)$$

Where $\epsilon$ represents the terms on the right-hand side of equation (1), $f_{111to100}(n)$ and $f_{111to100}(n)$ represents the transformation of ½<111> loop to <100> loop and <100> loop to ½<111> loop, respectively, during cascade overlap. Similarly, differential equations are also employed for <100> loops and voids.

For the present study, the cascade overlap rate is calculated based on approximation in terms of experimental flux and volume of the cascade-induced molten domains in the damaged region. In the previous work conducted by Backer et al. [56] using numerical models, cascade simulations are carried out to formulate a relationship which shows that there is a linear relationship between the total volume of molten domain ($V_{mol}$) and the PKA ($E_{PKA}$) energy, $V_{mol} = a*E_{PKA}$ and a =

$V_{fr}/E_{fr}$, $V_{fr}$ is volume of the molten domain at the cascade fragmentation energy ($E_{fr}$) which is 160 keV [56]. $V_{fr}$ is $10^6$ $A^{o3}$ at the threshold $E_{fr}$.

In order to obtain the $F_{ol}$, the neutron flux spectrum from the experiment (reported in Table 3) is employed in SPECTRA-PKA [58] code to analyse the PKA spectrum for range of PKA energy. The PKA spectrum is generated employing TENDL-2017 [59] pre-processed nuclear data which consists of recoil matrices for W shown in Figure 5 (a). Based on the linear relationship, $V_{mol} = a*E_{PKA}$, $V_{mol}$ is calculated using the $E_{PKA}$ and multiplied by the PKA spectrum for values above $E_{fr}$ to obtain average cascade overlap rate, $F_{ol}$, which is around $5 \times 10^{-7}/s$.

### 3.4 Comparison of CD model results with experimental data

Figure 5 (b-d) and Figure 7 depicts the total loop and void density distribution and their comparison between experiment and numerical results at irradiation temperatures of 400 °C and 600 °C. The results obtained from the numerical model is able to reproduce the trend of dislocation loops (½<111>+<100>) and voids from the experimental data in terms of the loop and void sizes (Figure 5 (b-d), Figure 7). However, an exact agreement is not achieved for loop and void density, which can be attributed to the various factors employed in the CD model. One reason is that treatment of trapped elements such interstitial impurities and substitutional transmuted element in the CD model. Since the concentration of impurities are less due to 99.7% pure W sample, the role of irradiation-induced transmuted elements can influence the numerical results. In order to study the actual migration of the transmuted elements and their interaction with the loops and voids, separate ODEs are required for transmuted elements, which need to be coupled in the CD model, which can be computationally expensive. Since the main objective of the current study is to analyse the long-term evolution of the loops and voids in W material, transmuted elements are only considered as traps. The other reason is the lack of data on the mixed and 'None' clusters and based on some assumptions, these clusters are integrated in ½<111>, <100> loops and c15 clusters as percentage fractions.

Based on experimental observations, mean size of loops (½<111>+<100>) of 5.3 nm and 2.3 nm, respectively, are observed. Moreover, percentage of <100> loops is measured as 53% and 21% at 400 °C and 600 °C, respectively. The remaining of the SIA clusters can be comprised of ½<111>, mixed, 'None' and c15 clusters and their measurements are not carried out. With regards to the temperature dependence, the total loop density distribution is found to be higher at 600 °C with respect to 400 °C. In the case of numerical results, loop sizes less than 19 nm are able to be obtained in the model at both irradiation temperatures (400 °C and 600 °C). <100> loop size up to 16 nm and ½<111> loop size less than 19 nm are predicted in the model as depicted in the Figure 5 (d). In particular, the total loop density at 600 °C is found to be higher than at 400 °C, which is consistent with the experimental observations. The density of ½<111> loops are larger than the <100> loops. At 600 °C, higher density distributions are obtained for ½<111> loops with respect to 400 °C. In the case of <100> loops, the scenario is quite different. At loop size greater than 7 nm, <100> loop density at 600 °C is lower than at 400 °C. This is due to the fact that <100> loops absorb smaller <111> loops to form larger <100> loops at 400 °C. Figure 6 (a-b) shows the evolution of the ½<111> and <100> loops, in terms of density and average size, as a function of irradiation dose. It is interesting to note the fact that at dpa less than 0.02, the density of ½<111> loop at 400 °C is higher than at 600 °C. While at dpa greater than 0.02, there is decrease in the ½<111> loop density at 400 °C with respect to 600 °C. This is due to the absorption of smaller ½<111> loops by <100> loop which reduces loop population of ½<111> loop and increases the <100> loop density. However, for <100> loops, the overall effect is that density at 600 °C is greater than at 400 °C. This is corroborated in the Figure 6 c in terms number of SIAs ($N_{SIAs}$) accommodated by ½<111> and <100> loops. The average size of the loops (½<111>+<100>) is predicted to be around 3.8- 4 nm in the CD model (Figure 6 b). Concerning c15 cluster evolution, it is worthwhile to note that at 600 °C, SIAs accommodated by loops (½<111>+<100>) are higher than in c15 clusters (Figure 6 d). However, at 400 °C and dpa greater than 0.02, more migrating SIAs are used for forming c15 clusters with consequent overall reduction of SIAs in the loops.

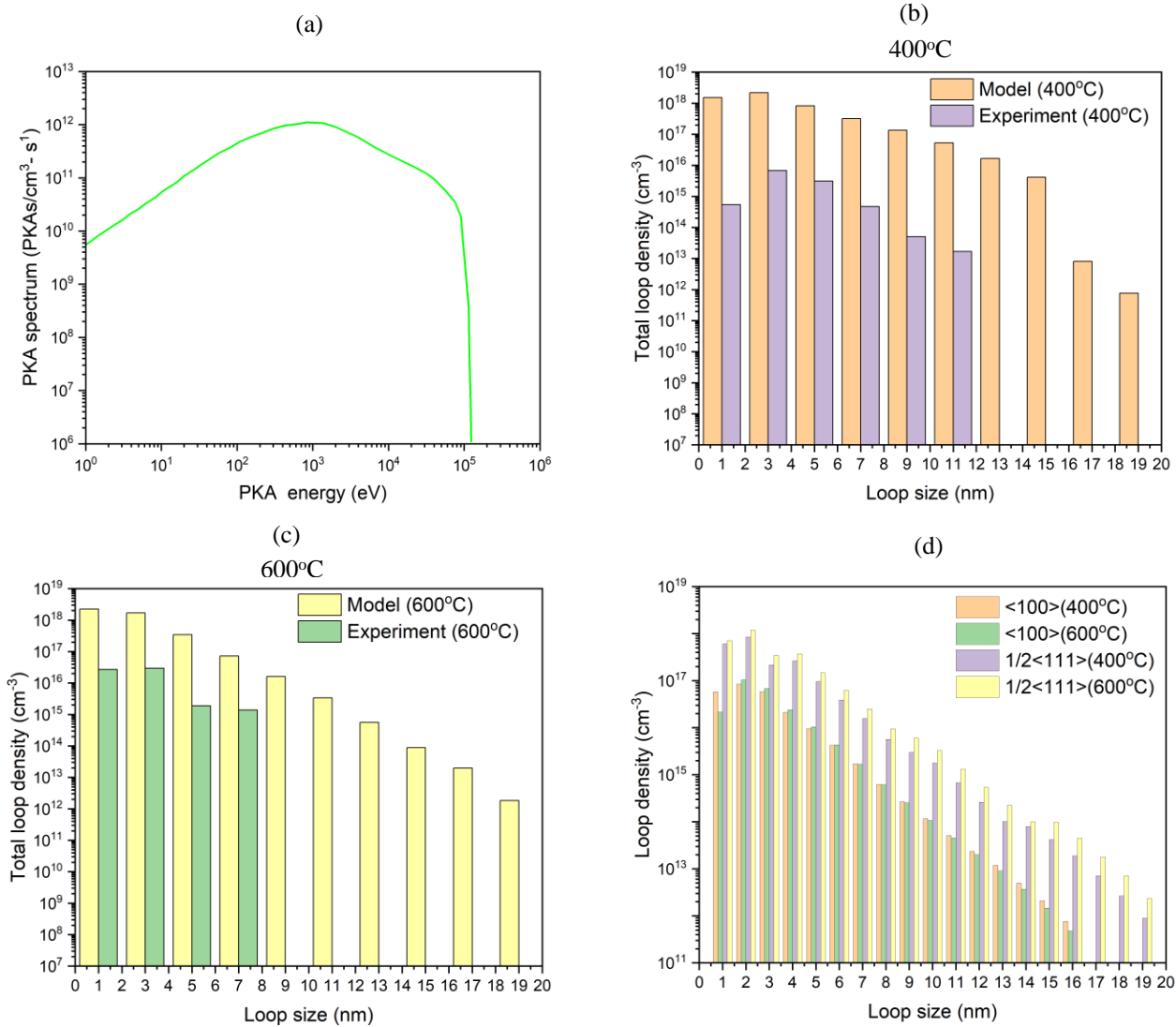

*Figure 5. (a) PKA spectrum from SPECTRA-PKA across the DCT specimen based on the experimental irradiation conditions. Total loop density vs loop size obtained from experimental and CD model at (b) 400 °C (c) 600 °C for 1.02 dpa. (d) 1/2<111> and <100> loop density from CD model at 400 °C and 600 °C for 1.02 dpa.*

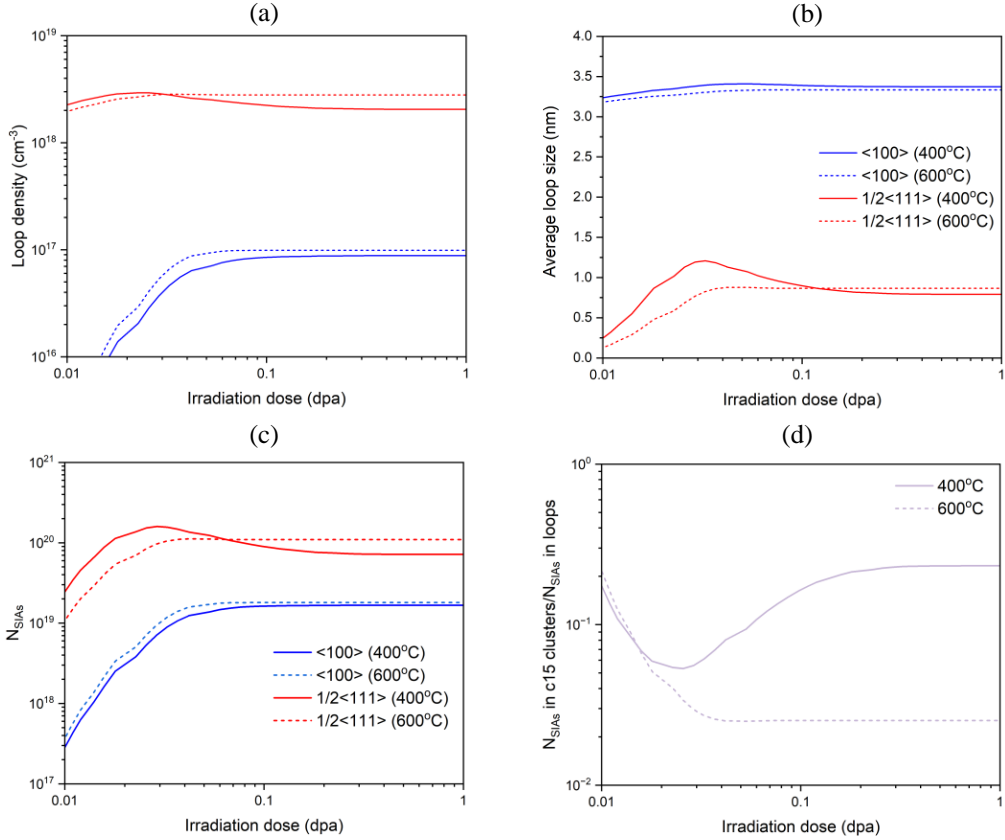

*Figure 6. CD model results as a function of irradiation dose (dpa) for a) ½<111>and <100> density b) Average loop size c) SIAs in ½<111> and <100> (in cm$^{-3}$) d) ratio between SIAs in c15 clusters and SIAs in ½<111> and <100>.*

In the case of voids, the model overpredicts the density distribution with respect to experimental observations (Figure 7). It should be noted that at 600 °C, the void size increases when compared to 400 °C. Moreover, the average size of voids is higher at 600 °C than at 400 °C as depicted in Figure 7 (bottom, right). At 400 °C, it is obvious that for void size greater than 1.5 nm, the difference in the void density between experimental and numerical results are higher Figure 7 (bottom, left). This may be due to parameters and assumptions implemented in the CD model as explained in the first paragraph of this subsection. The overall void density distribution at 600 °C is lower compared to 400 °C. This is mainly due to the fact that at higher temperatures, more SIAs coalesce with voids leading to reduction in the void density. In addition, the presence of trapped elements impedes the mobility of SIAs which consequently interacts with the voids and decreases

void density at smaller void sizes. In fact, the smaller mobile voids (<=4 vacancies) are absorbed by immobile voids leading to larger void size.

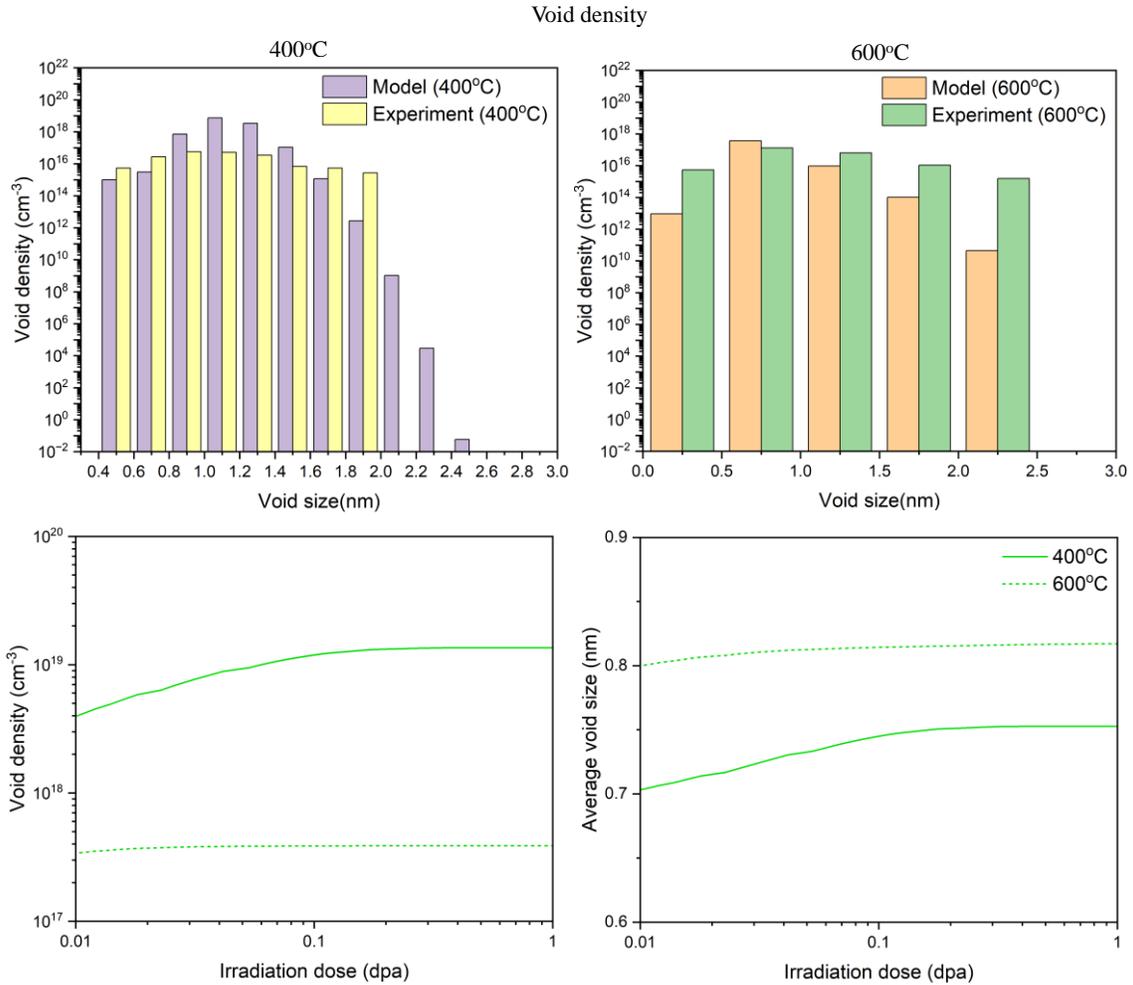

*Figure 7. Void density vs void size obtained from experimental and CD model at 400 °C and 600 °C for 1.02 dpa (top), CD model results as a function of irradiation dose (dpa) for void density (left) average void size (right) (bottom).*

*3.5    Influence on the Cascade overlap*

To analyse the overlap of the cascade with dislocation loops and voids and the consequent change in the overall defect population, simulations are performed, and a comparison is made between the cases with overlap event and no overlap event. In the cascade overlap, a fraction of ½<111> loops, <100> loops and void are transformed to <100> loops, <100> loops and voids, respectively, and it is modelled in this study by means of transformation coefficients as explained in the section 3.2. Figure 8 (a-c) shows the loop and void density distribution of different scenarios at temperatures 400 °C and 600 °C. It is obvious that in the event of cascade overlap at 400 °C, there is an increase in the ½<111> loop density distribution (1.24 x $10^{18}$ to 2.04 x $10^{18}$ $cm^{-3}$) while for <100> loops, there is a small decrease in the size upto 10 nm which increases slightly afterwards. The number of SIAs accommodated by ½<111> loops are higher than the <100> loops. An increase in the ½<111> loop size in noticed for cascade overlap events since the ½<111> loops absorb smaller <100> loops and therefore, there is an increase in the ½<111> loop size. However, at higher temperature of 600 °C, the scenario is different. There is only smaller increase in the ½<111> total loop density for cascade overlap case (2.73 x $10^{18}$ to 2.9 x $10^{18}$ $cm^{-3}$). Moreover, the increase in the density and size of <100> loops are more pronounced. This can be due to the fact that at higher temperatures along with the transformation of ½<111> to <100> loop during cascade overlap, the <100> loops subsume smaller ½<111> loops to elevate its density and size. There is an overall decrease in the void density at 400 °C and 600 °C (Figure 8 (c)). This is attributed to the fact that the migrating SIAs formed during the cascade overlap event interact and recombine with the voids, thus reducing the void density.

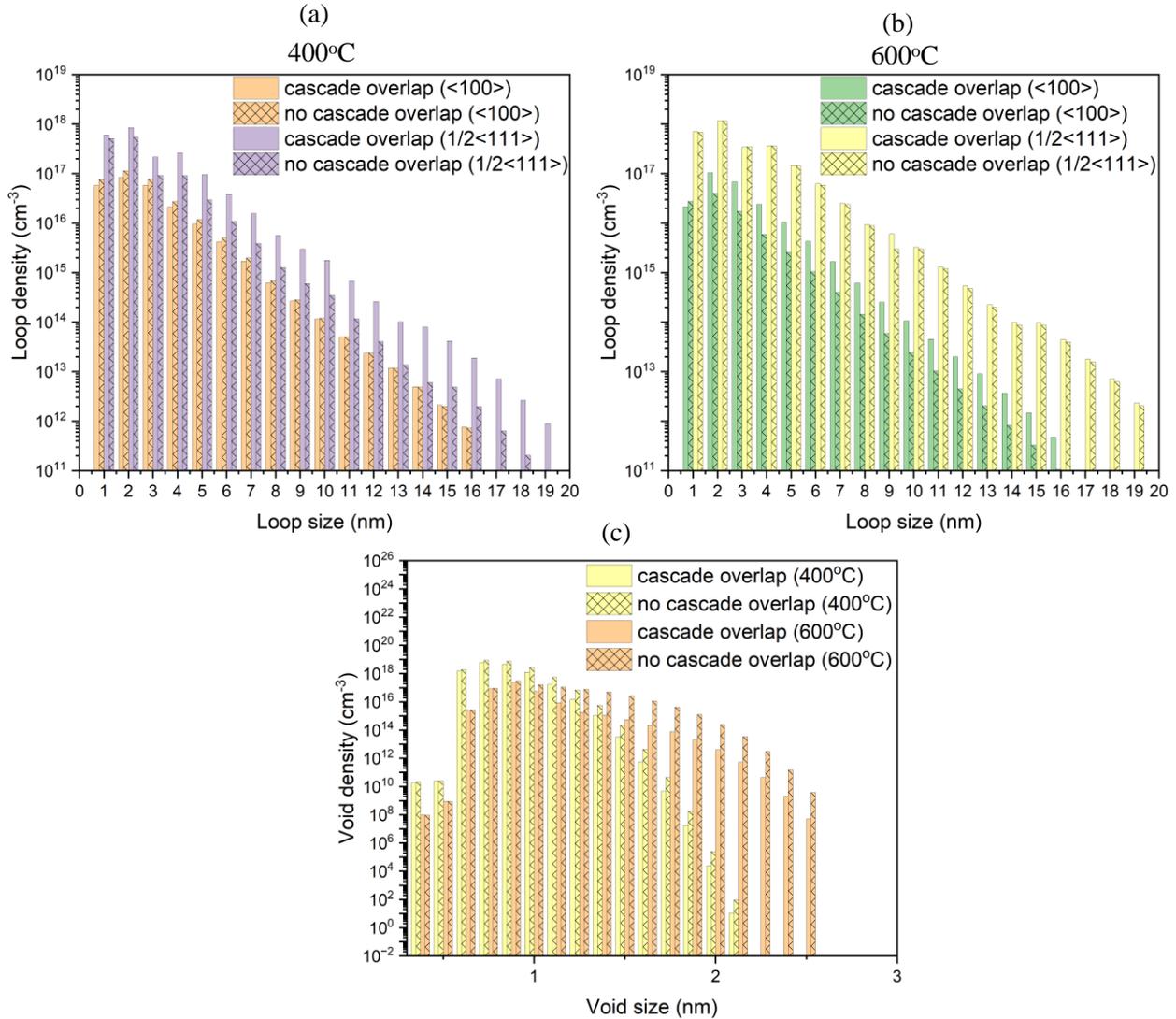

*Figure 8. Influence of Loop density vs loop size from CD model due to cascade overlap at (a) 400 °C and (b) 600 °C for 1.02 dpa (c) Influence of void density vs void size from CD model due to cascade overlap at 400 °C and 600 °C for 1.02 dpa.*

*3.6    Influence of loops and voids due to trapped elements.*

Figure 9 depict the ½<111> loop, <100> loop and void density distribution for irradiation temperatures 400 °C and 600 °C. The results demonstrate the comparison between the presence (with traps) and absence (without traps) of the trapped elements. It is obvious from the in the presence of interstitial impurity and substitutional transmuted elements, higher density and size distribution of ½<111> loops are observed, which are trapped by the interstitial and transmuted elements. In the case without the inclusion of the trap elements, the ½<111> loop density distribution is less. Moreover, ½<111> loop size of less than 5 nm is noticed without trap elements. Due to the high mobility of ½<111> loops in W, part of the ½<111> loops during their long-range movement are absorbed at sink (grain boundaries, dislocation line) and most of the remaining ½<111 loops are absorbed by <100> loops. In the case of <100> loops, their density and size in the absence of trap elements are lower than with respect to the case with the presence of trap elements. However, a higher density and size of <100> loops are obtained in the presence of trap elements when compared to ½<111> loops. This is due to the fact that, as mentioned before, most of the ½<111> loops are subsumed by <100> loops. It is interesting to note the fact that in the presence of trap elements at 400 °C, most of the mobile ½<111> loops are impeded by the presence of traps. Therefore, their coalescence with the voids is reduced and a higher density and size (<=2.3 nm) of voids is noticed. In the case of without trap elements, high mobile ½<111> loops have the possibility to coalesce with voids, which leads to the reduction in the density and size of voids (<=1.5 nm) when compared to the case with interaction of ½<111> loops with traps. A similar scenario is also noticed at the temperature of 600 °C.

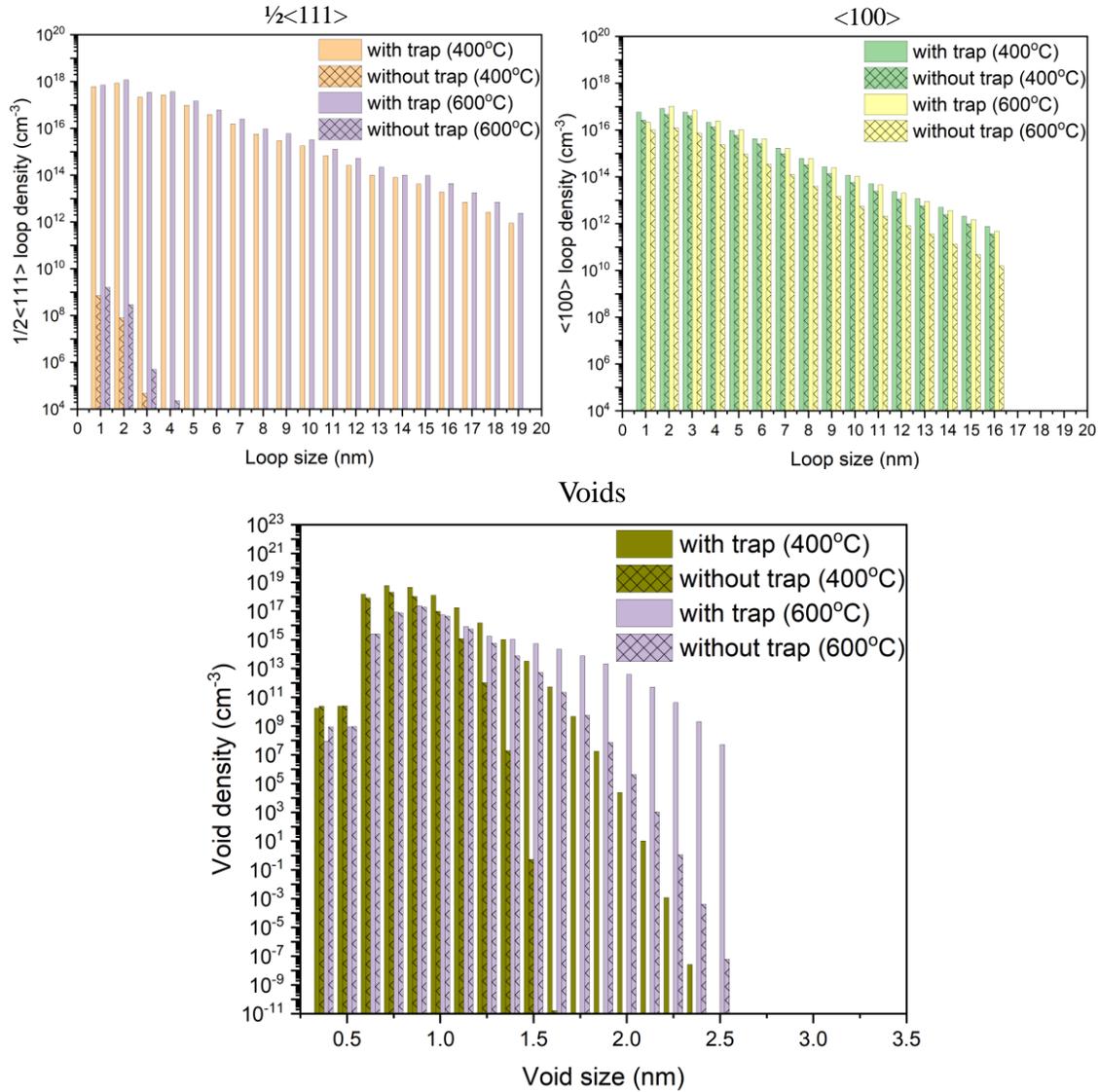

*Figure 9. Influence due to trapped elements on 1/2<111> and <100> loop density from CD model for 1.02 dpa (top), influence due to trapped elements on void density from CD model for 1.02 dpa (bottom).*

Figure 10 (a-b) depicts the evolution of ½<111> and <100> loops represented in terms of number of SIAs ($N_{SIAs}$) in loops (½<111> and <100>) as a function of irradiation dose and temperature. It is worthwhile to note the fact that the SIAs accommodated by ½<111> loops in the presence of traps are much higher than the case without traps. However, for <100> loops, some of ½<111> loops are absorbed at sink without trap. Therefore, the ½<111> loops subsumed by <100> loops are lower with smaller number of SIAs employed for <100> loop formation. The SIAs, accommodated in the evolution of c15 cluster, are greatly influenced by the presence of trap elements. The number of SIAs for the c15 clusters (Figure 10 d) is lower in the presence of trap elements since most of the SIAs are used in the ½<111> and <100> loop formation at temperature of 400 °C and 600 °C. When there are no traps involved, the highly mobile ½<111> loops interact with voids, which decreases the void density with respect to the presence of traps (Figure 10 c). Moreover, with the increase in the irradiation dose, elevation in the void density is observed (Figure 10 c).

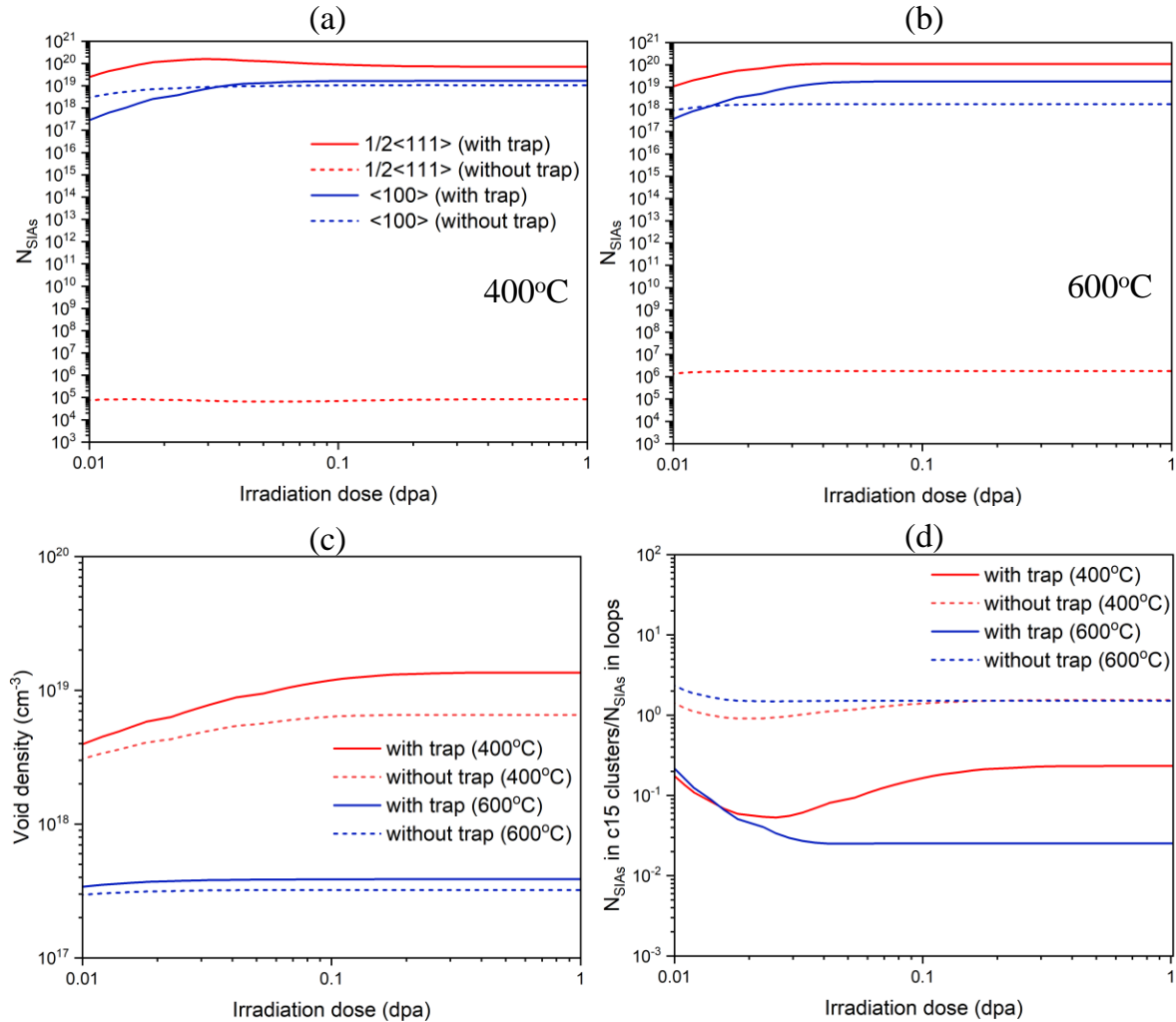

*Figure 10. CD model results as a function of irradiation dose (dpa) for SIAs in ½<111> and <100> at a) 400 °C b) 600 °C (c) void density d) ratio between SIAs in c15 clusters and SIAs in ½<111> and <100>.*

## 3.7 Influence of loop reaction

In order to understand the influence of ½<111> loop and <100> loop interaction on the overall defect population in W material, simulations are performed for reaction 'on' and reaction 'off' conditions. Reaction 'on' corresponds to the reaction between ½<111> and <100> loop while reaction 'off' does not take into account the interaction between ½<111> and <100> loops. Figure 11 depicts the loop and void density distribution in reaction 'on' and reaction 'off' conditions. Figure 12 (a-b) depicts the $N_{SIAs}$ in the ½<111> and <100> loop for irradiation temperatures 400 °C and 600 °C. Moreover, ratio between SIAs in ½<111> and SIAs in <100> loops are calculated and depicted in Figure 12 b. To analyse the density of c15 clusters, ratio between SIAs in c15 clusters and SIAs in ½<111> + <100> loops are also evaluated as shown in Figure 12 d.

For irradiation temperatures 400 °C and 600 °C, the density distribution of ½<111> loops is higher in reaction 'off' condition when compared to the reaction 'on' condition as Figure 11. In particular, the higher density of ½<111>loops are more pronounced at a loop size greater than 2 nm. At irradiation temperature of 400 °C, <100> loop density is lower in reaction 'off' conditions for loop size between 3 nm and 8 nm. This is due to the fact that the reaction between ½<111> loops and <100> loops form larger <100> loops at size greater than 8 nm which consequently decreases the density of <100> loops for size between 3 nm and 8 nm. In general, for reaction 'on' condition, the interaction between ½<111> and <100> loops promote the growth of <100> loops 400 °C and 600 °C (Figure 12 (a-b)). At higher temperature of 600ºC, the reaction 'off' decreases <100> loop density when compared to temperature at 400 °C. However, for reaction 'on' condition, smaller <100> loops are absorbed by ½<111> loops which consequently reduces overall population <100> loops. The $N_{SIAs}$ accommodated by the ½<111> and <100> loops are evaluated based on the ratio between SIAs in ½<111> and SIAs in <100> loops (Figure 12 (c)). In reaction 'on' condition, the reaction between ½<111> and <100> loops promote the <100> loop growth by accommodating higher $N_{SIAs}$ with consequent decrease in the ratio at dpa greater than 0.08 dpa. However, for reaction 'off' conditions, at 400 °C, <100> loop grows by absorbing the smaller <100> loops and migrating SIAs during in-cascade clustering and cascade overlap, which increase the loop size with consequent decrease in the ratio of the $N_{SIAs}$. In case of 600 °C, loop growth of <100> loop is impeded due to the unavailability of additional $N_{SIAs}$ for the loop growth, which lead to the increase in the ratio as depicted in figure. Since there are no transfer of SIAs between ½<111> and <100> loops in reaction 'off' condition, more SIAs are accommodated in the c15 clusters with respect to reaction 'on' condition (Figure 12 d). Moreover, at higher temperatures, more SIAs are accommodated by loops than c15 cluster. The density of voids is influenced by the SIAs involved in the formation of the loops. Since there are no transfer of SIAs between ½<111> and <100> loops, most of the SIAs will either be absorbed by c15 clusters or voids. It is quite obvious from that lower density of voids are observed in reaction 'off' condition at the temperature 400 °C and 600 °C at size greater than 1 nm since the SIAs from loops coalesce with the vacancies.

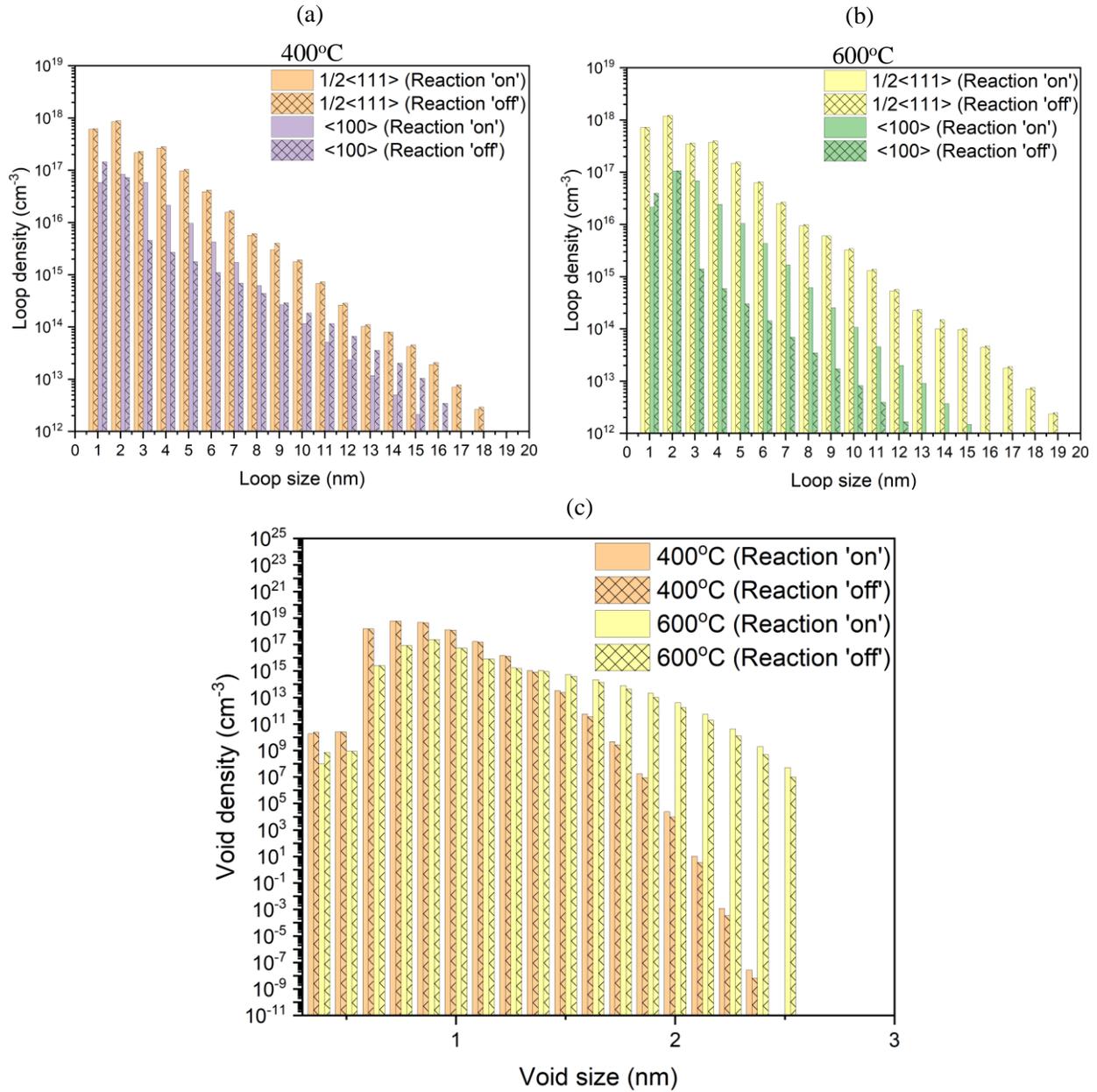

*Figure 11. Influence of loop density vs loop size from CD model due to reaction conditions between ½<111> and <100> loops at (a) 400 °C (b) 600 °C for 1.02 dpa. (c) Influence of void density vs void size from CD model due to reaction conditions between ½<111> and <100> loops at 400 °C and 600 °C for 1.02 dpa.*

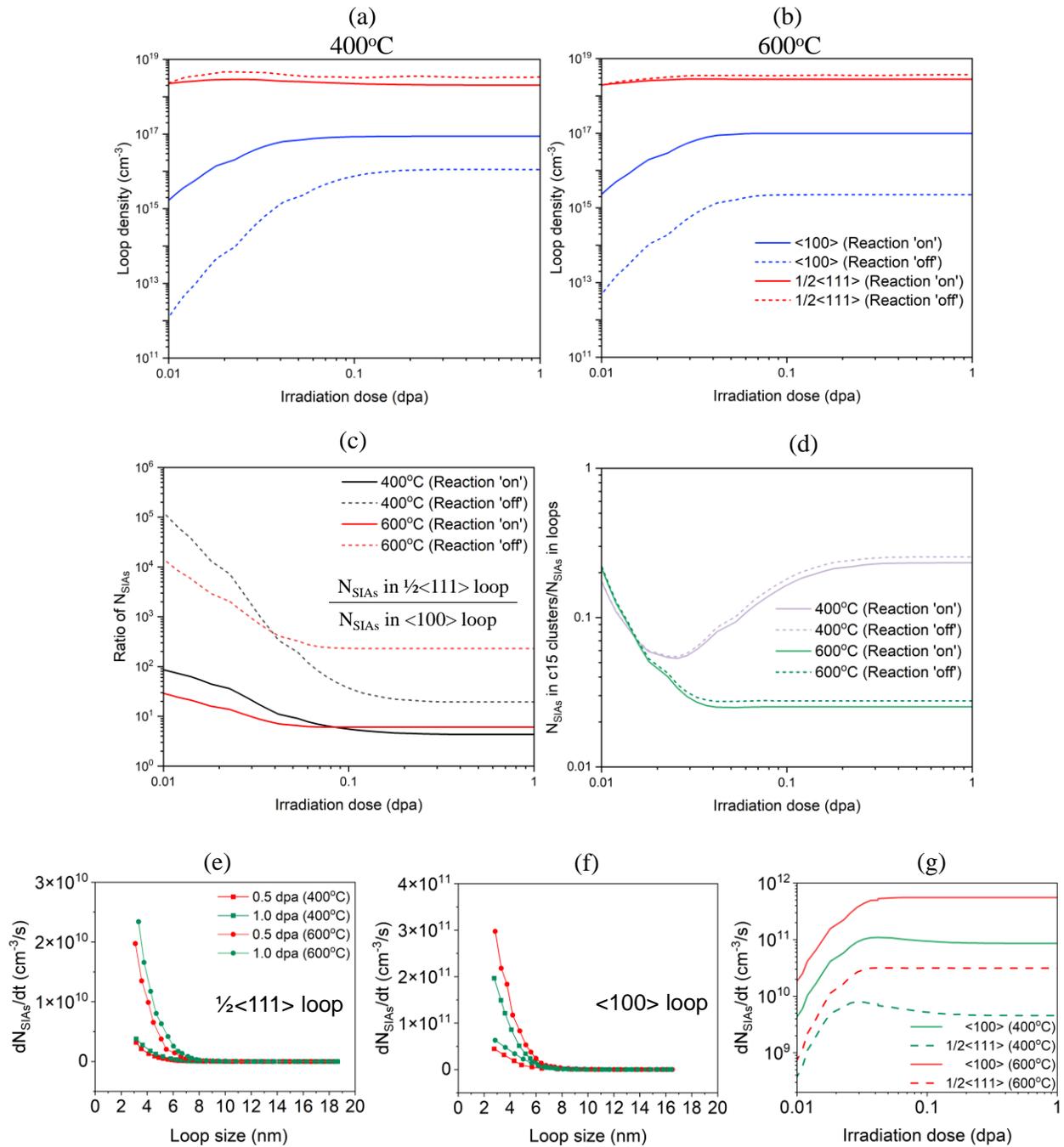

*Figure 12. CD model results as a function of irradiation dose (dpa) for ½<111> and <100> density at a) 400 °C b) 600 °C c) ratio between SIAs in ½<111> and <100> d) ratio between SIAs in c15 clusters and SIAs in ½<111> and <100>.*

## 4  Discussions

There are several mechanisms involved in the formation of <100> loops in irradiated W. One such mechanism in the formation of <100> loops is the interaction between <111> loops of identical size. However, such kind of reaction is rare. Since the formation energy of ½<111> loops are lower than <100> loops, the probability for the formation of <100> loops is lower with respect to the formation of the ½<111> loops [14]. In order to increase the probability of interaction between ½<111> loop variants, factors such as extended irradiation time, elevation in irradiation temperature are required [14]. The model has considered a probability of ∼3/8 for the formation of <100> loops resulting from the reaction between ½<111> loops of same size but the proportion of <100> loops formed through this interaction is very less in the present study. Therefore, in the current model, the main source for the formation of the <100> loops are the generation directly in irradiation cascades and full cascade overlap with pre-existing defects. The density of ½<111> loops are found to be higher than <100> loops at various irradiation temperatures. The formation of 1/2<111> loops and <100> loops by means of c15 cluster collapse is not considered in the study since the irradiation temperature employed in the experimental and numerical study is not sufficient for the collapse of c15 cluster as explained in the previous work of Liu et al. [38,39]. However, there are only few studies available in the literature regarding the stability of c15 clusters, which can influence the long-term evolution of ½<111> and <100> density and size.

The temporal evolution of the defect features is influenced by the interaction with the trap elements (impurity atoms + transmuted elements). Due to these interactions, the size and density of loops and voids have influence on the material behaviour as explained in the section 3.5. Since the sample used in this study is a bulk material, the surface effects are not important. Therefore, the highly mobile ½<111> loops can be sunk at grain boundaries or be incorporated into dislocation lines. In particular, with the inclusion of the trap elements, the model can capture higher concentration of ½<111> loops. Such kind of scenario is highlighted in the experimental work of Klimenkov et al. [16]. Moreover, a detailed study is required to include the inclusion of Re atoms in the CD model to understand a long-range interaction of Re atoms with loops and voids.

In order to increase the rate of interaction between ½<111> and <100> loops, the influence of mean free path is relevant, which depends on the loop population of ½<111> and <100> loops. To understand the evolution of the loops at extended irradiation fluence, there is a need to gain some insights into the mutual transfer of SIAs due to the reaction between ½<111> and <100> loops.

The overall loop population actually depends on the reaction between ½<111> and <100> loops, which fosters the loop growth via the absorption or coalescence reactions. Due to this, there is a decrease in the total population of the loops. In fact, the change in the loop population results in the SIA transfer between the loops and influence the SIAs in the loops with respect to its initial microstructure. The SIA transfer rate between the ½<111 and <100> loops is studied based on similar approach carried out by Gao et al. [37] and the expressions for SIA transfer rate for ½<111> and <100 loops is calculated as follows in equation 19 and equation 20:

$$Tran_{111}(j) = \sum_{j>i} k_i^j C_{111}(j) C_{100}(i).i - \sum_{j<i'} k_{i'}^j C_{111}(j) C_{100}(i').j \qquad (19)$$

$$Tran_{100}(j) = \sum_{j>i} k_i^j C_{100}(j) C_{111}(i).i - \sum_{j<i'} k_{i'}^j C_{100}(j) C_{111}(i').j \qquad (20)$$

The SIAs received by ½<111> loops through absorbing the smaller loop size of <100> type is depicted in the $k_i^j C_{111}(j) C_{100}(i).i$ term, while $k_{i'}^j C_{111}(j) C_{100}(i').j$ term depicts the <100> loop of larger size absorbing the smaller loop size ½<111> type. A positive value for $Tran_{111}(j)$ depicts the gain of SIAs by ½<111> of size j from <100> loops through '½<111>+ <100>' reaction. In the case of negative $Tran_{111}(j)$, there is a loss of SIAs from ½<111>loop. The terms in equation 19 and 20 have similar depiction for <100 > loops.

Figure 12 (e-f) depicts the SIA transfer rate of ½<111 and <100 loops. It is worthwhile to note in general both ½<111>and <100 loops gain SIAs due to their reaction due to the positive Tran$_{111}$ and Tran$_{100}$, which corroborates the fact that the SIA transfer promotes the loop growth of ½<111> and <100>. The transfer rate of <100> loop is higher with respect to ½<111> loop. In particular, for the smaller loops sizes less than 6 nm, the SIA transfer rate is higher when compared to the larger loop size for ½<111> and <100> loop types. <100> loop of sizes less than 6 nm achieves larger gain of SIAs due to the higher production of smaller ½<111 loops as a result of SIAs aggregation and in-cascade clustering. The lesser density of the higher sized loops (½<111>, <100>) are due to the low transfer rate of SIAs. Moreover, with the increase in the irradiation dose and temperature, the SIA transfer rate is increased for ½<111> loop. At temperature of 600 °C, the reduction in SIA gain of <100> loops promote an increase in gain in SIAs for <111> loops. An assessment of the loop growth is also carried out by means of total SIA transfer rates (Figure 12

g) by summing up the transfer rates of ½<111> and <100> loop types and plotted as a function of irradiation dose. At 400 °C, the transfer rate of <100> loops are higher with respect to ½<111> loops. Moreover, there is a higher reduction in the SIA transfer rate after 0.02 dpa for ½ <111> loop when compared to <100> loops. However, at 600 °C, the transfer rate is constant for ½<111> and <100> loops after 0.02 dpa with higher transfer rates for <100> loop.

## 5    Conclusions

Neutron irradiation-induced defects pose the main challenge in the design and development of the in-vessel components in a fusion reactor due to alteration of the material properties. Moreover, the studies addressing behaviour of W based components in fusion relevant conditions are still ongoing. In order to understand the microstructural evolution of irradiation-induced defects in W, the present study included an integrated experimental-numerical analysis on neutron irradiated W. Regarding experimental analysis, TEM characterisation is carried out for the W irradiated up to 1 dpa at 400 °C and 600 °C. CD method, developed in this study, is employed for investigating the long-term evolution of ½<111> loops, <100> loops, voids and c15 clusters by assigning the irradiation conditions implemented for experimental analysis. The numerical results obtained from the CD method are compared with the experimental data. The following statements highlight the findings obtained in this work:

1) Based on the comparison between experimental and computational results, the CD model provides insights on the evolution of dislocation loops and voids at two irradiation temperatures (400 and 600 °C) and could aid in understanding the evolution of microstructural defect features in fusion relevant conditions. In particular, temperature dependence of the populations of voids and loops is clearly observed in both experimental and numerical results. The discrepancy between experimental and numerical results in terms of the size and density of void and loops can be due to the treatment of mixed and 'None' clusters which are modelled as a part of ½<111> and <100> loops. Moreover, separate concentration equation has to be employed for the transmuted atoms in the CD model which will be the following work based on the current study.

2) Cascade overlap study is performed by incorporating the transformation coefficients in CD

model to evaluate the full overlap event with pre-existing SIA clusters and the consequent change in the long-term evolution of the defects is analysed. In particular, the study has shown that the cascade overlap serves as the additional source for the formation of <100> loops.

3) Trap mediated diffusion of ½<111> loops and mobile vacancies (<=4) are modelled in the CD model by including interstitial impurity and substitutional atoms as traps in the CD model. The influence of traps on the concentration of dislocation loops and voids is observed with respect to the case without the inclusion of traps in the CD model.

4) The CD results have shown that the reaction between ½<111> loops and <100> loops promote the SIA transfer between the loops and influence the loop growth of both loop types. In fact, the SIA transfer rate of <100> loops found to be higher than the ½<111> loops.

Based on the previous studies employing the CD model, a key tenet of the current work is the assessment of the long-term evolution of ½<111> loop, <100> loop and voids, especially for larger defect sizes. The hope is that the parameters and results obtained from this integrated experimental-numerical approach can be useful for providing information for upper scale models like dislocation dynamics, finite element model for investigating the engineering properties and thermo-mechanical behaviour of fusion relevant components.

**Acknowledgment:**

This work is supported by the Helmholtz Association of German Research Centers (HGF) and has been carried out within the framework of Nuclear Fusion Programme at Karlsruhe Institute of Technology (KIT), Germany. Authors would like to thank all the members of the Fusion Materials Laboratory (FML) of the Karlsruhe Institute of Technology (KIT), Germany for their help in handling irradiated materials.

This work has been carried out within the framework of the EUROfusion Consortium, funded by the European Union via the Euratom Research and Training Programme (Grant Agreement No 101052200 ─ EUROfusion). Views and opinions expressed are however those of the author(s) only and do not necessarily reflect those of the European Union or the European Commission. Neither the European Union nor the European Commission can be held responsible for them.

**CRediT authorship contribution statement**

**Salahudeen Mohamed:** Conceptualization, Methodology Software, Formal analysis, Investigation, Writing – original draft, Visualization. **Qian Yuan:** Resources, Investigation, Writing – review & editing. **Dimitri Litvinov:** Resources, Writing – review & editing. **Jie Gao:** Investigation, Software, Writing – review & editing. **Ermile Gaganidze:** Investigation, Resources, Writing – review & editing. **Dmitry Terentyev:** Investigation, Writing – review & editing. **Hans-Christian Schneider:** Investigation, Writing – review & editing. **Jarir Aktaa:** Investigation, Resources, Writing – review & editing.

**Declaration of Competing Interest**

The authors declare that they have no known competing financial interests or personal relationships that could have appeared to influence the work reported in this paper.